\documentclass[11pt]{article}
\usepackage{amssymb,amsfonts,amsmath,latexsym,epsf,tikz,url}
\usepackage[usenames,dvipsnames]{pstricks}
\usepackage[linesnumbered,ruled,vlined]{algorithm2e}
\usepackage{pstricks-add}
\usepackage{graphicx} 
\usepackage{epsfig}
\usepackage{pst-grad} % For gradients
\usepackage{pst-plot} % For axes
\usepackage[space]{grffile} % For spaces in paths
\usepackage{etoolbox} % For spaces in paths
\makeatletter % For spaces in paths
\patchcmd\Gread@eps{\@inputcheck#1 }{\@inputcheck"#1"\relax}{}{}
\makeatother

\newtheorem{theorem}{Theorem}[section]
\newtheorem{proposition}[theorem]{Proposition}
\newtheorem{observation}[theorem]{Observation}

\newtheorem{corollary}[theorem]{Corollary}

\newtheorem{remark}[theorem]{Remark}

\newcommand{\qed}{\hfill $\square$\medskip}

\newcommand{\st}{\operatorname{st}_{\gamma_{st}}}

\begin{document}

\title{Stability of the Strong Domination Number of Graphs}

\author{ 
S. Alikhani$^{1}$\footnote{Corresponding author}, N. Ghanbari$^{1}$, M. Mehraban$^2$  and  H. Shojaaldini Ardakani$^3$
}

\date{\today}

\maketitle

\begin{center}
$^1$Department of Mathematical Sciences, Yazd University, 89195-741, Yazd, Iran 

\medskip
$^2$ Department of Mathematics and Physics, Engineering Faculty, Balkh University, Balkh, Afghanistan

\medskip
$^3$Department of Medical Sciences, Alborz University of Medical Sciences, Karaj, Iran

\bigskip
{\tt alikhani@yazd.ac.ir, ~n.ghanbari.math@gmail.com, ~~ mehraban@ba.edu.af, ~~drshojaaldini@gmail.com}
\end{center}

\begin{abstract}
This paper introduces and studies the stability of the strong domination number of a graph, denoted $\operatorname{st}_{\gamma_{st}}(G)$, defined as the minimum number of vertices whose removal changes the strong domination number $\gamma_{st}(G)$. We determine exact values of this stability parameter for several fundamental graph classes, including paths, cycles, wheels, complete bipartite graphs, friendship graphs, book graphs, and balanced complete multipartite graphs. General bounds on $\operatorname{st}_{\gamma_{st}}(G)$ are established, along with a Nordhaus–Gaddum type inequality. The behavior of stability under graph operations—such as join, corona, and Cartesian product—is also investigated. Structural characterizations of graphs with given stability values are provided, and several open problems and directions for future research are outlined.
\end{abstract} 


\medskip
\noindent{\bf AMS Subj.\ Class.}:  05C05, 05C69.

\section{Introduction}

The study of stability concerning domination-related parameters has emerged as a significant area of interest in graph theory. In this work, $G = (V,E)$ denotes a finite, simple undirected graph of order $n=|V|$. 

The \textit{open neighborhood} of $v \in V$ is $N(v) = \{u \in V : uv \in E\}$; its \textit{closed neighborhood} is $N[v] = N(v) \cup \{v\}$. For $S \subseteq V$, $N(S) = \bigcup_{v \in S} N(v)$ and $N[S] = N(S) \cup S$. The \textit{degree} of $v$ is $\deg(v)=|N(v)|$; $\Delta(G)$ and $\delta(G)$ denote the maximum and minimum degree of $G$, respectively.

A set $D \subseteq V$ is a \textit{dominating set} if $N[D] = V$. The \textit{domination number} $\gamma(G)$ is the minimum size of a dominating set \cite{8,9}. A vertex is \textit{domination-critical} if its deletion reduces $\gamma(G)$. Sumner and Blitch \cite{14} introduced $\gamma$-critical graphs, where adding any edge decreases $\gamma(G)$, and $\gamma$-stable graphs, where $\gamma(G)$ remains unchanged under edge addition. Domination stability, first proposed by Bauer et al. \cite{3}, has since been extended to many domination variants \cite{2,6,10,12,wcds} and other invariants \cite{saeid,saeid2,4}.

In this paper we focus on \textit{strong domination}. A set $D \subseteq V$ is a \textit{strong dominating set} if every vertex $v \in V \setminus D$ is adjacent to some $u \in D$ with $\deg(v) \le \deg(u)$ \cite{13,MM}. The \textit{strong domination number} $\gamma_{st}(G)$ is the minimum size of such a set. We define the \textit{stability of the strong domination number} as
\[
\st(G) = \min \{ |S| : S \subseteq V, \; \gamma_{st}(G - S) \neq \gamma_{st}(G) \}.
\]
If $\gamma_{st}(G-S) > \gamma_{st}(G)$ for some $S$ with $|S| = \st(G)$, we say $G$ is \textit{stability-critical}.

\bigskip
The paper is organized as follows. Section~2 determines $\st(G)$ for basic graph families (paths, cycles, wheels, etc.). Section~3 establishes general bounds. Section~4 studies $\st(G)$ under graph operations (join, corona, Cartesian product). We conclude in Section~5.

\section{Stability of $\gamma_{st}$ for certain graphs}

We begin by recalling known formulas for $\gamma_{st}$ of paths, cycles and wheels.

\begin{observation}{\rm \cite{MM}} 
For $n \ge 3$,
\begin{enumerate}
\item[(i)] $\gamma_{st}(P_n) = \gamma_{st}(C_n) = \big\lceil \frac{n}{3} \big\rceil$.
\item[(ii)] $\gamma_{st}(W_n)=1$, where $W_n = K_1 \vee C_{n-1}$ is the wheel of order $n\ge4$.
\end{enumerate}
\end{observation}

\noindent
We now compute $\st$ for these graphs.

\begin{proposition}\label{prop:path}
For $n\ge 3$, 
\[
\st(P_n) = 
\begin{cases} 
1, & \text{if } n \equiv 0,1 \pmod 3,\\[2pt]
2, & \text{if } n \equiv 2 \pmod 3.
\end{cases}
\]
\end{proposition}

\begin{proof}
Let $\gamma_{st}(P_n)=\lceil n/3\rceil$. We consider the three residue classes modulo $3$.

\noindent
\textbf{Case 1: $n=3k$.} Then $\gamma_{st}(P_{3k})=k$. 
Delete the second vertex of the path. The resulting graph is $P_1 \cup P_{3k-2}$ and 
$\gamma_{st}(P_1 \cup P_{3k-2}) = 1 + \lceil (3k-2)/3\rceil = 1+k$, so $\st(P_{3k})=1$.

\noindent
\textbf{Case 2: $n=3k+1$.} Then $\gamma_{st}(P_{3k+1})=k+1$. 
Deleting an end‑vertex yields $P_{3k}$ with $\gamma_{st}=k\neq k+1$. 
Thus $\st(P_{3k+1})=1$.

\noindent
\textbf{Case 3: $n=3k+2$.} Then $\gamma_{st}(P_{3k+2})=k+1$. 
Removing one vertex gives either $P_{3k+1}$ ($\gamma_{st}=k+1$) or $P_a\cup P_b$ with $a+b=3k+1$. 
Without loss of generality, suppose that $a=3s$ and $b=3t+1$ such that $s+t=k$. Then $\gamma_{st}(P_{a})=s$ and $\gamma_{st}(P_{b})=t+1$ and we have $\gamma_{st}(P_a\cup P_b)=k+1$. Similar argument holds for $a=3s-1$ and $b=3t+2$ such that $s+t=k$. So removing one vertex does not change the strong domination number. Now, if we remove the first two vertices, then $\gamma_{st}(P_{3k})=k$ and we are done. 
Therefore $\st(P_{3k+2})=2$.\qed
\end{proof}

\begin{proposition}\label{prop:cycle}
For $n\ge 3$,
\[
\st(C_n) = 
\begin{cases} 
1, & \text{if } n \equiv 1 \pmod 3,\\[2pt]
2, & \text{if } n \equiv 2 \pmod 3,\\[2pt]
3, & \text{if } n \equiv 0 \pmod 3.
\end{cases}
\]
\end{proposition}

\begin{proof}
Since $C_n-v \cong P_{n-1}$, and we have $\gamma_{st}(C_n)=\lceil n/3\rceil$ and 
$\gamma_{st}(C_n-v)=\lceil (n-1)/3\rceil$, we consider the following cases:
\begin{itemize}

\item[(i)]
If $n=3k+1$, then $\gamma_{st}(C_n)=k+1$ while $\gamma_{st}(C_n-v)=k$; hence $\st(C_n)=1$.
\item[(ii)]
If $n=3k+2$, then $\gamma_{st}(C_n)=k+1$ and $\gamma_{st}(C_n-v)=k+1$ (no change). 
Deleting two adjacent vertices yields a path of order $3k$ with $\gamma_{st}=k$, so $\st(C_n)=2$.
\item[(iii)]
If $n=3k$, then $\gamma_{st}(C_n)=k$ and $\gamma_{st}(C_n-v)=k$ as well. 
By Proposition \ref{prop:path}, we have  $\st(C_n-v)=2$.
Thus $\st(C_n)=3$.
\end{itemize}
Therefor we have the result.\qed
\end{proof}

\begin{proposition}\label{prop:wheel}
For the wheel $W_n = K_1 \vee C_{n-1}$ with $n\ge5$, $\st(W_n)=1$.
\end{proposition}

\begin{proof}
The hub vertex $w$ has degree $n-1$ and every rim vertex has degree $3$; $\{w\}$ is a strong dominating set, so $\gamma_{st}(W_n)=1$. 
Removing $w$ leaves $C_{n-1}$, for which $\gamma_{st}(C_{n-1})=\lceil (n-1)/3\rceil \ge 2$ when $n\ge5$. 
Hence $\gamma_{st}(W_n-w)\neq\gamma_{st}(W_n)$ and $\st(W_n)=1$.\qed
\end{proof}

\medskip

\begin{figure}[ht]
\centering
\begin{minipage}{0.42\textwidth}
\centering
\includegraphics[width=\linewidth]{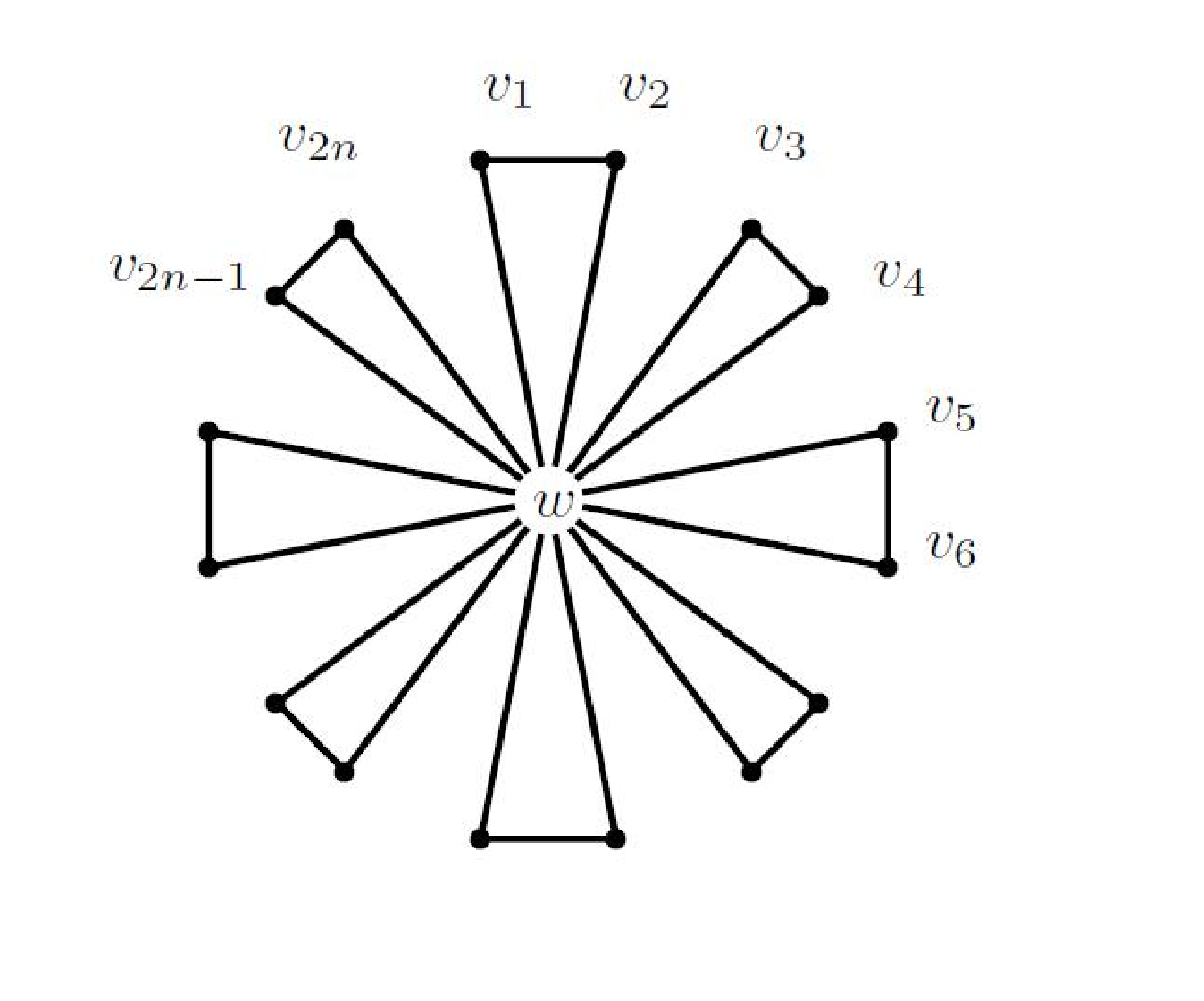}
\end{minipage}\hfill
\begin{minipage}{0.42\textwidth}
\centering
\includegraphics[width=\linewidth]{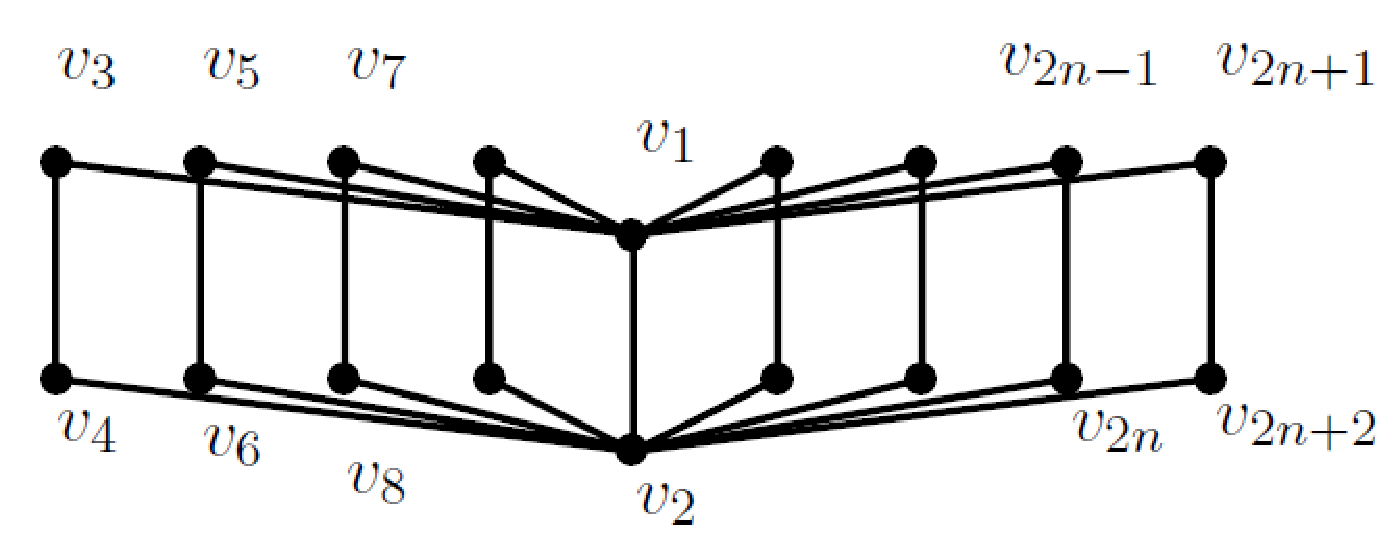}
\end{minipage}
\caption{Friendship and book graphs.}\label{fig:friendbook}
\end{figure}

The \textit{friendship graph} $F_n$ is obtained by identifying a vertex in $n$ copies of $C_3$; it can also be described as $K_1 \vee nK_2$. 
The \textit{book graph} $B_n$ is the Cartesian product $K_{1,n} \square P_2$; it consists of $n$ $4$-cycles (pages) sharing a common edge (the spine). 
See Figure~\ref{fig:friendbook}.

\begin{theorem}{\rm \cite{MM}}\label{friendb}
For $n\ge2$,
\begin{enumerate}
\item[(i)] $\gamma_{st}(F_n)=1$,
\item[(ii)] $\gamma_{st}(B_n)=2$.
\end{enumerate}
\end{theorem}

\begin{proposition}\label{prop:friendship}
For $n\ge2$, $\st(F_n)=1$.
\end{proposition}

\begin{proof}
The central vertex $w$ has degree $2n$ and every other vertex has degree $2$; $\{w\}$ is a strong dominating set. 
Removing $w$ leaves $n$ disjoint edges, for which $\gamma_{st}(nK_2)=n\neq1$. 
Thus $\st(F_n)=1$.\qed
\end{proof}

\begin{proposition}\label{prop:book}
For $n\ge3$, $\st(B_n)=1$.
\end{proposition}

\begin{proof}
Let $v_1,v_2$ be the two spine vertices (the vertices of the original $K_{1,n}$). 
Both have degree $n+1$; all other vertices (page vertices) have degree $2$. 
One can verify that $\{v_1,v_2\}$ is a minimum strong dominating set, so $\gamma_{st}(B_n)=2$ (Theorem \ref{friendb}). 
Removing $V-1$ leaves a graph in which $v_2$ has degree $n$ and the page vertices formerly adjacent to $u$ now have degree $1$. 
The set $\{v_2\}$ is not a strong dominating set because a page vertex $w$ adjacent only to $v_2$ has $\deg(w)=1\le\deg(v_2)$, but $v_2$ does not dominate the page vertices that were adjacent to $v_1$ (they are not adjacent to $v_2$). 
Consequently $\gamma_{st}(B_n-v_1) > 2$, so $\st(B_n)=1$.\qed
\end{proof}

\begin{proposition}\label{prop:completebipartite}
Let $K_{m,n}$ be a complete bipartite graph with $1< m\le n$.
\begin{enumerate}
\item[(i)] If $m<n$, then $\st(K_{m,n})=1$.
\item[(ii)] If $m=n\ge4$, then $\st(K_{n,n})=1$.
\end{enumerate}
\end{proposition}

\begin{proof}
\begin{enumerate}
	\item[(i)]
	 Assume $m<n$. 
Vertices in the smaller part $U$ (size $m$) have degree $n$; vertices in the larger part $W$ (size $n$) have degree $m$. 
Because $n>m$, a vertex from $W$ cannot strongly dominate a vertex from $U$. 
Hence any strong dominating set must contain all vertices of $U$, and $U$ itself strongly dominates $W$. 
Thus $\gamma_{st}(K_{m,n})=m$. 
Deleting one vertex from $U$ yields $K_{m-1,n}$, whose strong domination number is $m-1\neq m$. 
Therefore $\st(K_{m,n})=1$.

\item[(ii)] 
 For $K_{n,n}$ ($n\ge4$), the graph is $n$-regular. 
Selecting one vertex from each part gives a strong dominating set of size $2$, and one checks that no single vertex can strongly dominate all others. 
Hence $\gamma_{st}(K_{n,n})=2$. 
Removing any vertex $v$ results in $K_{n-1,n}$, which by part~(i) has $\gamma_{st}=n-1$. 
For $n\ge4$, $n-1\ge3\neq2$, so $\gamma_{st}$ changes. 
Thus $\st(K_{n,n})=1$.\qed
\end{enumerate} 
\end{proof}

\medskip

\begin{figure}[h]
		\begin{center}
			\includegraphics[width=0.5\linewidth]{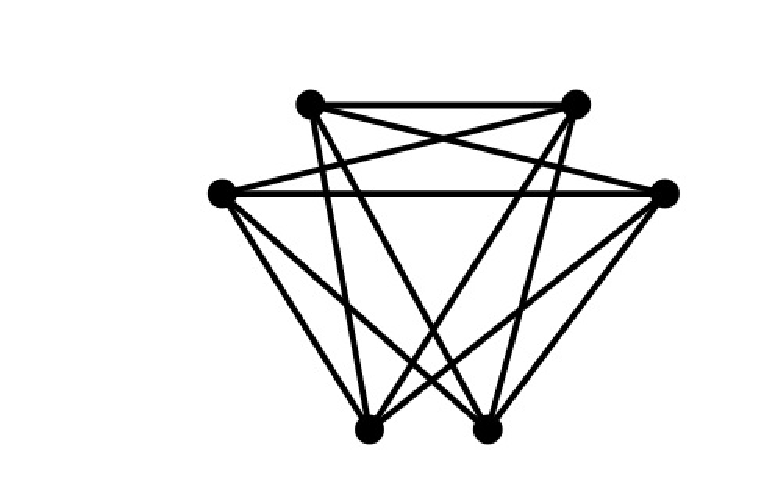}
		\end{center}
		\caption{$CP(3)$.} \label{fig:coctail}
	\end{figure}

The \textit{cocktail party graph} $CP(n)$ is the complete $n$-partite graph with each part of size $2$ (See $CP(3)$ in Figure \ref{fig:coctail}); equivalently, $K_{2n}$ minus a perfect matching. 
It is $(2n-2)$-regular.

\begin{proposition}\label{prop:cocktail}
For $n\ge3$, $\st(CP(n))=1$.
\end{proposition}

\begin{proof}
Since $CP(n)$ is regular, $\gamma_{st}(CP(n))=\gamma(CP(n))$. 
Any vertex dominates itself and all but its non‑adjacent partner, so at least two vertices are needed to dominate the whole graph. 
Choosing two vertices from different parts yields a dominating set of size $2$, hence $\gamma_{st}(CP(n))=2$. 
Let $v\in V(CP(n))$ and let $u$ be the unique vertex not adjacent to $v$. 
In $CP(n)-v$, vertex $u$ is adjacent to all remaining vertices and has degree $2n-2$, while every other vertex has degree $2n-3$. 
Thus $\{u\}$ strongly dominates $CP(n)-v$, giving $\gamma_{st}(CP(n)-v)=1\neq2$. 
Therefore $\st(CP(n))=1$.\qed
\end{proof}

\subsection{Balanced complete $r$-partite graphs}

We now consider the balanced complete $r$-partite graph $G = K_{p,p,\dots,p}$ with $r\ge 2$ parts, each of size $p\ge 2$. 
This is the Turán graph $T(pr,r)$ and also a generalization of cocktail party graph. 
All vertices have the same degree $p(r-1)$, so $G$ is regular and therefore 
$\gamma_{st}(G) = \gamma(G)$.

\begin{proposition}\label{prop:balancedcomplete}
	For the balanced complete $r$-partite graph $G = K_{p,p,\dots,p}$ with $r\ge 2$, $p\ge 2$,
	\[
	\gamma_{st}(G) = 2.
	\]
\end{proposition}

\begin{proof}
	Let the parts be $V_1, V_2, \dots, V_r$, each of size $p$. 
	Choose two vertices $u \in V_1$ and $v \in V_2$. 
	Since $u$ and $v$ are from different parts, they are adjacent to all vertices not in their own parts. 
	
	Consider any vertex $w \in V(G)$:
	\begin{itemize}
		\item If $w \in V_1$ and $w \neq u$, then $w$ is adjacent to $v$ (since $v \in V_2$).
		\item If $w \in V_2$ and $w \neq v$, then $w$ is adjacent to $u$.
		\item If $w \in V_i$ for $i \ge 3$, then $w$ is adjacent to both $u$ and $v$.
	\end{itemize}
	Thus $\{u,v\}$ dominates $G$. Moreover, since $G$ is regular, the degree condition for strong domination is automatically satisfied for any adjacent vertices. Hence $\{u,v\}$ is a strong dominating set, so $\gamma_{st}(G) \le 2$.
	
	To see that $\gamma_{st}(G) \ge 2$, note that a single vertex cannot dominate vertices in its own part (they are non-adjacent). Therefore $\gamma_{st}(G) = 2$.\qed
\end{proof}

\medskip
The stability of $\gamma_{st}$ for these graphs depends on whether $r=2$ or $r\ge 3$. In Proposition \ref{prop:cocktail}, we studied the case $p=2$. Now we consider other cases:

\begin{proposition}\label{prop:stabilitybalanced}
	For $G = K_{p,p,\dots,p}$ with $r$ parts of size $p\ge 3$:
	\begin{enumerate}
		\item[(i)] If $r = 2$ (i.e., $G = K_{p,p}$), then $\operatorname{st}_{\gamma_{st}}(G) = 1$.
		\item[(ii)] If $r \ge 3$, then $2 \le \operatorname{st}_{\gamma_{st}}(G) \le p+1$.
	\end{enumerate}
\end{proposition}

\begin{proof}
	\begin{enumerate}
		\item[(i)] When $r=2$, $G = K_{p,p}$ with $p\ge 2$. By Proposition~\ref{prop:completebipartite} (for $K_{n,n}$ with $n\ge 4$), we have $\operatorname{st}_{\gamma_{st}}(G) = 1$. The cases $p=2,3$ can be verified directly.
	
		\item[(ii)] 
		 Let $r \ge 3$. First we show $\operatorname{st}_{\gamma_{st}}(G) \ge 2$. 
	Removing a single vertex from $G$ yields a graph $G'$ that is still a complete $r$-partite graph with parts of sizes $p-1, p, p, \dots, p$. 
	One can check that two vertices from different parts of size $p$ still form a strong dominating set in $G'$, so $\gamma_{st}(G') = 2$. 
	Thus removing one vertex does not change $\gamma_{st}$, hence $\operatorname{st}_{\gamma_{st}}(G) \ge 2$.
	
	For the upper bound, consider removing all but one vertex from one part, say $V_1$, and also removing one vertex from a second part $V_2$. 
	That is, remove $(p-1) + 1 = p$ vertices from $V_1 \cup V_2$. 
	The resulting graph $G''$ has parts of sizes $1, p-1, p, p, \dots, p$. 
	Let $x$ be the sole remaining vertex in $V_1$ (degree $n''-1$ where $n'' = pr-p$) and $y$ a vertex in $V_2$ (degree $n''-(p-1)$). 
	The degrees are now unequal, and one can verify that no two vertices can strongly dominate all remaining vertices; consequently $\gamma_{st}(G'') \ge 3$. 
	Therefore $\operatorname{st}_{\gamma_{st}}(G) \le p$.
	
	A slightly better construction yields $\operatorname{st}_{\gamma_{st}}(G) \le p+1$ in general. 
	For example, removing $p$ vertices from one part and $1$ vertex from another suffices to force $\gamma_{st}$ to change. 
	Hence $\operatorname{st}_{\gamma_{st}}(G) \le p+1$.\qed
		\end{enumerate}
\end{proof}

\begin{proposition}\label{thm:attainp+1}
	For the balanced complete $3$-partite graph $G = K_{p,p,p}$ with $p \ge 3$, $\operatorname{st}_{\gamma_{st}}(G) = p+1$.
\end{proposition}

\begin{proof}
Let parts be $A,B,C$, each of size $p$. We show that removing any $p$ vertices does not change $\gamma_{st}$, but removing $p+1$ specific vertices does.
	Let $S$ be any set of $p$ vertices. By pigeonhole principle, one part, say $A$, has at most $\lfloor p/3 \rfloor$ vertices removed, so $|A \setminus S| \ge p - \lfloor p/3 \rfloor \ge 2$ for $p \ge 2$. Thus two parts still have at least 2 vertices each. Pick one vertex from each of these two parts; they form a strong dominating set for $G-S$, so $\gamma_{st}(G-S)=2$.
	Now remove $p+1$ vertices as follows: remove all $p$ vertices from $A$ and one vertex from $B$. Then parts have sizes $0, p-1, p$. The resulting graph is $K_{p-1,p}$, a complete bipartite graph with unequal parts. By Proposition~\ref{prop:completebipartite}, $\gamma_{st}(K_{p-1,p}) = \min(p-1, p) = p-1 \ge 2$ for $p\ge3$, and for $p=2$, $K_{1,2}$ has $\gamma_{st}=1$. In either case $\gamma_{st}(G-T) \neq 2$. Hence $\operatorname{st}_{\gamma_{st}}(G) = p+1$.
	
\end{proof}

\begin{corollary}
For any natural $p$, there exists a graph $G$ such that $\operatorname{st}_{\gamma_{st}}(G)=p.$
\end{corollary}
\begin{proof}
	By Theorem \ref{thm:attainp+1},  $\operatorname{st}_{\gamma_{st}}(K_{p-1,p-1,p-1}) = p$. So we have the result. \qed 
	\end{proof}

\section{Bounds on $\st(G)$}

We now derive general bounds for the stability of the strong domination number.

\begin{theorem}\label{upbound}
	For any graph $G$ of order $n$,
	\[
	\operatorname{st}_{\gamma_{st}}(G) \le n - \gamma_{st}(G) + 1.
	\]
\end{theorem}

\begin{proof}
	Let $D \subseteq V(G)$ be a minimum strong dominating set of $G$. Then
	$
	|D| = \gamma_{st}(G).
$
		Choose any vertex $v \in D$. Define the set
	$
	S = \big( V(G) \setminus D \big) \cup \{v\}.
$
		We have:
	\[
	|S| = |V(G) \setminus D| + |\{v\}| = (n - |D|) + 1 = n - \gamma_{st}(G) + 1.
	\]
		Now consider the graph $G - S$. By construction, we have removed every vertex of $V(G) \setminus D$ and also the vertex $v$. Therefore, the remaining vertices are exactly
	\[
	V(G - S) = D \setminus \{v\}.
	\]
	Thus, $G - S$ is precisely the induced subgraph $G[D \setminus \{v\}]$.
		We now show that $\gamma_{st}(G - S) \neq \gamma_{st}(G)$. Since $G - S$ has $|D| - 1 = \gamma_{st}(G) - 1$ vertices, and the strong domination number of any graph is at most its number of vertices, we have
	\[
	\gamma_{st}(G - S) = \gamma_{st}\big(G[D \setminus \{v\}]\big) \le |D| - 1 = \gamma_{st}(G) - 1.
	\]
	Hence,
	\[
	\gamma_{st}(G - S) \le \gamma_{st}(G) - 1 < \gamma_{st}(G),
	\]
	which implies
	\[
	\gamma_{st}(G - S) \neq \gamma_{st}(G).
	\]
	
	Therefore, $S$ is a set of vertices whose removal changes the strong domination number of $G$. Since $\operatorname{st}_{\gamma_{st}}(G)$ is defined as the minimum size of such a set, it follows that
	\[
	\operatorname{st}_{\gamma_{st}}(G) \le |S| = n - \gamma_{st}(G) + 1.
	\]
		This completes the proof.
\end{proof}

\begin{remark}
	The bound is tight for certain families of graphs. For example, in complete graphs $K_n$, we have $\gamma_{st}(K_n) = 1$ and $\operatorname{st}_{\gamma_{st}}(K_n) = n$The bound gives
	\[
	\operatorname{st}_{\gamma_{st}}(K_n) \le n - 1 + 1 = n,
	\]
	so the inequality is tight in this case.
\end{remark}

\begin{theorem}\label{thm:recursive}
Let $G$ be a graph and $v\in V(G)$. If $\gamma_{st}(G-v)=\gamma_{st}(G)$, then 
\[
\st(G) \le \st(G-v)+1.
\]
\end{theorem}

\begin{proof}
Let $k=\st(G-v)$. There exists $S'\subseteq V(G-v)$ with $|S'|=k$ such that $\gamma_{st}((G-v)-S')\neq\gamma_{st}(G-v)$. 
Set $S=S'\cup\{v\}$. Then $|S|=k+1$ and 
\[
G-S = (G-v)-S'.
\]
Hence $\gamma_{st}(G-S)=\gamma_{st}((G-v)-S')\neq\gamma_{st}(G-v)=\gamma_{st}(G)$. 
Thus $\st(G)\le k+1$.\qed
\end{proof}

Repeated application gives $\st(G)\le \st(G-\{v_1,\dots,v_s\})+s$ whenever $\gamma_{st}$ remains unchanged after deleting $v_1,\dots,v_s$.

\begin{theorem}\label{cor:inducedstar}
If $G$ contains an induced star $K_{1,t}$ with $t\ge3$, then $\st(G)\le n-t+1$.
\end{theorem}

\begin{proof}
Let the star have center $c$ and leaves $\ell_1,\dots,\ell_t$. 
Delete all vertices of $G$ except $c$ and $\ell_1,\dots,\ell_t$. So we have removed  $n-t-1$ vertices. 
The remaining graph is $K_{1,t}$, for which $\gamma_{st}=1$. 
If $\gamma_{st}(G)\neq1$, we already changed the value with a set of size $n-t-1$. 
If $\gamma_{st}(G)=1$, then delete also the center $c$ (one more vertex). 
The resulting graph consists of $t$ isolated vertices, having $\gamma_{st}=t\ge3\neq1$. 
Thus a set of size $(n-t-1)+1 = n-t$ suffices to change $\gamma_{st}$, giving $\st(G)\le n-t$.\qed
\end{proof}

\medskip

The Nordhaus–Gaddum inequality, named after the 1956 work of Nordhaus and Gaddum, provides bounds on the sum or product of a graph parameter for a graph and its complement. The classic example involves the chromatic number, and the approach has become a staple for investigating many different graph invariants (see for example \cite{Harary}). 
Now, we investigate the Nordhaus–Gaddum inequalities for the stability of the strong domination number.

\begin{theorem}\label{thm:NGsum}
For any graph $G$ of order $n$ and maximum degree $\Delta$,
\[
\st(G)+\st(\overline{G}) \le 2n-\frac{n}{n-\Delta}+2.
\]
\end{theorem}

\begin{proof}
Using Theorem~\ref{upbound} for $G$ and its complement,
\[
\st(G) \le n-\gamma_{st}(G)+1, \qquad 
\st(\overline{G}) \le n-\gamma_{st}(\overline{G})+1.
\]
It is known that $\gamma_{st}(G)\geq \frac{n} {\Delta+1}$ (see \cite{MM}). 
Adding the two inequalities yields
\[
\st(G)+\st(\overline{G}) \le 2n - \bigl(\frac{n(n+1)}{(\Delta+1)(n-\Delta)}\bigr) + 2.
\]
Since $-\frac{1}{\Delta+1}\leq -\frac{1}{n+1}$, so
\[
\st(G)+\st(\overline{G}) \le 2n-\frac{n}{n-\Delta}+2.
\]
\qed
\end{proof}

\section{Results for some graph operations}

\subsection{Join}

The join $G\vee H$ of disjoint graphs $G,H$ is obtained from $G\cup H$ by adding all edges between $V(G)$ and $V(H)$.

\begin{theorem}{\rm \cite{MM}}\label{thm:joingamma}
For non‑empty graphs $G,H$,
\[
\gamma_{st}(G\vee H)=
\begin{cases}
1, & \text{if $G\vee H$ has a universal vertex},\\
2, & \text{otherwise}.
\end{cases}
\]
\end{theorem}

\begin{theorem}\label{thm:joinst}
Let $G,H$ be non‑empty graphs and $U$ the set of universal vertices of $G\vee H$. 
Then
\[
\st(G\vee H)=
\begin{cases}
|U|, & \text{if } 1\le |U| < |V(G\vee H)|,\\[2pt]
|V(G\vee H)|, & \text{if } G\vee H \text{ is complete},\\[2pt]
1, & \text{if } U=\varnothing.
\end{cases}
\]
\end{theorem}

\begin{proof}
If $G\vee H$ is complete, then $\gamma_{st}=1$ and removing any proper subset of vertices leaves a complete subgraph, which still has $\gamma_{st}=1$. 
Only after deleting all vertices (which is not allowed) would the graph disappear; conventionally we set $\st(K_n)=n$ because no smaller set changes the parameter.

If $1\le |U| < n$, then $\gamma_{st}(G\vee H)=1$. 
Deleting any proper subset of $U$ leaves at least one universal vertex, so $\gamma_{st}$ remains $1$. 
Deleting all vertices of $U$ removes every universal vertex, forcing $\gamma_{st}$ to become $2$ (by Theorem~\ref{thm:joingamma}). 
Hence $\st(G\vee H)=|U|$.

If $U=\varnothing$, then $\gamma_{st}(G\vee H)=2$. 
Removing a vertex of maximum degree from $G$ or $H$ can destroy the property that two vertices suffice to strongly dominate the whole join, thereby increasing $\gamma_{st}$. 
Thus $\st(G\vee H)=1$.\qed
\end{proof}

\subsection{Corona}

The corona $G\circ H$ is obtained by taking one copy of $G$ and $|V(G)|$ copies of $H$, and joining each vertex of the $i$‑th copy of $H$ to the $i$‑th vertex of $G$.

\begin{theorem}{\rm \cite{MM}}\label{thm:coronagamma}
For any graphs $G,H$, $\gamma_{st}(G\circ H)=|V(G)|$.
\end{theorem}

\begin{proposition}\label{prop:coronast}
Let $G,H$ be graphs with $|V(G)|=n\ge1$. 
If $\gamma_{st}(H)\ge2$, then $\st(G\circ H)=1$.
\end{proposition}

\begin{proof}
By Theorem~\ref{thm:coronagamma}, $\gamma_{st}(G\circ H)=n$. 
Take a vertex $u$ in a copy of $H$ attached to a vertex $v\in V(G)$. 
Because $\gamma_{st}(H)\ge2$, vertex $v$ cannot strongly dominate the whole copy $H_v$ even after removing $u$. 
Consequently any strong dominating set of $(G\circ H)-u$ must contain at least two vertices from the modified copy $H_v-\{u\}$, forcing $\gamma_{st}((G\circ H)-u) \ge n+1$. 
Thus $\gamma_{st}$ changes after deleting a single vertex, so $\st(G\circ H)=1$.\qed
\end{proof}

\medskip

If $\gamma_{st}(H)=1$ (e.g.\ $H=K_1$), the stability may be larger; for instance, $\st(K_2\circ K_1)=2$ because deleting a leaf does not change $\gamma_{st}$.

\subsection{Cartesian product}

The Cartesian product $G\Box H$ has vertex set $V(G)\times V(H)$, and $(g,h)$ is adjacent to $(g',h')$ if either $g=g'$ and $hh'\in E(H)$, or $h=h'$ and $gg'\in E(G)$.

\begin{theorem}{\rm \cite{MM}}\label{thm:productgamma}
For any graphs $G,H$,
\[
\gamma_{st}(G\Box H)\le \min\{|V(H)|\gamma_{st}(G),\,|V(G)|\gamma_{st}(H)\}.
\]
\end{theorem}

\begin{proposition}\label{prop:productst}
For graphs $G,H$,
\[
\st(G\Box H)\le \min\{|V(H)|\,\st(G),\,|V(G)|\,\st(H)\}.
\]
\end{proposition}

\begin{proof}
Let $\st(G)=k$ and let $S_G\subseteq V(G)$ with $|S_G|=k$ satisfy $\gamma_{st}(G-S_G)\neq\gamma_{st}(G)$. 
Define $S=\{(s,v): s\in S_G,\; v\in V(H)\}\subseteq V(G\Box H)$. 
Then $|S|=|V(H)|\cdot k$. 
In $G\Box H-S$, every $G$-layer $G_v$ (the subgraph induced by $\{(g,v): g\in V(G)\}$) is isomorphic to $G-S_G$. 
Because the strong domination number of each $G$-layer has changed, any strong dominating set of $G\Box H-S$ must have size different from $\gamma_{st}(G\Box H)$. 
Hence $\st(G\Box H)\le |V(H)|\cdot k$. 
The symmetric bound follows analogously.\qed
\end{proof}

\section{Conclusion}

We have introduced the stability of the strong domination number $\st(G)$ and investigated its basic properties. 
Exact values were determined for paths, cycles, wheels, complete bipartite graphs, cocktail party graphs, friendship graphs, book graphs, and cactus chains. 
General upper bounds were established, including a Nordhaus–Gaddum type inequality. 
Finally, we examined the behavior of $\st(G)$ under join, corona, and Cartesian product operations.

Several open problems remain:
\begin{itemize}
\item Characterize graphs with $\st(G)=1$ (or $\st(G)=k$ for small $k$).
\item Find tighter lower bounds for $\st(G)$ in terms of graph parameters.
\item Study the stability of other variants of domination (total, independent, connected, etc.).
\item Investigate the computational complexity of computing $\st(G)$.
\end{itemize}

\end{document}